\documentclass[a4paper,12pt]{amsart}

\newtheorem{theorem}{Theorem}[section]
\newtheorem{lemma}[theorem]{Lemma}
\theoremstyle{definition}

\newtheorem{definition}[theorem]{Definition}

\theoremstyle{remark}

\normalsize
\numberwithin{equation}{section}

\begin{document}

\title{Universal series by trigonometric system in weighted $L^1_{\mu}$ spaces}
\author{S. A. Episkoposian}

\address{Department of Physics, Chair of Higher Mathematics, Yerevan State University, Yerevan, Al. Manoogian st.1,  375049, Armenia}
\thanks{The author was supported in part by Grant- 01-000 from the Government of Armenia .}


\date{}


\begin{abstract}
In this paper we consider the question of existence  of
trigonometric series universal in weighted $L^1_{\mu}[0,2\pi]$
spaces with respect to rearrangements and in usual sense.
\end{abstract}

\maketitle

\section{Introduction}
\par\par\bigskip
\par\par\bigskip

Let $X$ be a Banach space.
\begin{definition}
 A series
\begin{equation}
\sum_{k=1}^\infty f_k,\ \ f_k \in X
\end{equation}
is said to be universal in $X$ with respect to rearrangements, if
for any $f \in X$ the members of (1.1) can be rearranged so that
the obtained series $\displaystyle \sum_{k=1}^\infty
f_{\sigma(k)}$ converges to $f$ by norm of $X$.
\end{definition}

\begin{definition}
The series (1.1) is said to be universal (in $X$) in the usual
sense, if for any $f \in X$ there exists a growing sequence of
natural numbers $n_k$  such that the sequence of partial sums with
numbers $n_k$ of the series (1.1) converges to $f$ by norm of $X$.
\end{definition}

\begin{definition}
The series (1.1) is said to be universal (in $X$) concerning
partial series, if for any $f \in X$ it is possible to choose a
partial series $\displaystyle {\sum_{k=1}^\infty f_{n_k}}$ from
(1.1), which converges to the $f$ by norm of $X$.
\end{definition}

Note, that many papers are devoted (see [1]- [9]) to the question
on existence of variouse types of universal series in the sense of
convergence {\sl almost everywhere and on a measure}.

The first usual universal in the sense of convergence almost
everywhere trigonometric series were constructed by D.E.Menshov
[1] and V.Ya.Kozlov [2]. The series of the form
\begin{equation}
\frac {1}{2}+\sum_{k=1}^\infty a_k cos{kx}+b_k sin{kx}
\end{equation}
was constructed just by them such that for any measurable on
$[0,2\pi]$  function $f(x)$ there exists the growing sequence of
natural numbers $n_k$  such that the series (1.2) having the
sequence of partial sums with numbers $n_k$   converges to $f(x)$
almost everywhere on $[0,2\pi]$. (Note here, that in this result,
when $f(x)\in{L^1_{[0,2\pi]}} $, it is impossible to replace
convergence almost everywhere by convergence in the metric
${L^1_{[0,2\pi]}}$).

 This result was distributed by A.A.Talalian on arbitrary orthonormal complete systems (see [3]).
He also established (see [4]), that if $\{\phi_n(x)\}_{n=1}^\infty
$  - the normalized basis of space ${L^p_{[0,1]}},p>1 $, then
there exists a series of the form
\begin{equation}
\sum_{k=1}^\infty{a_k\phi_k(x)},\ \ a_k \to 0.
\end{equation}
which has property: for any measurable function $f(x)$ the members
of series (1.3) can be rearranged so that the again received
series converge on a measure on [0,1] to $f(x)$.

W. Orlicz [5] observed the fact that there exist functional series
that are universal with respect to rearrangements in the sense of
a.e. convergence in the class of a.e. finite measurable functions.

It is also useful to note that even Riemann proved that every
convergent numerical series which is not absolutely convergent is
universal with respect to rearrangements in the class of all real
numbers.

Let $\mu(x)$ be a measurable on $[0,2\pi]$ function with $
0<\mu(x) \le1, x\in[0,2\pi]$ and let $L_\mu^1[0,2\pi]$ be a space
of mesurable functions $f(x),\ \ x\in [0,2\pi]$ with
$$ \int_0^{2\pi} |f(x)| \mu(x) dx<\infty.$$

M.G.Grigorian constructed a series of the form (see [6]),
\[
\sum_{k=-\infty}^\infty{C_ke^{ikx} \ \  with \ \
\sum_{k=-\infty}^\infty \left | {C_k} \right|^q <\infty},\ \
\forall q>2
\]
which is universal in $L_\mu^1[0,2\pi]$ concerning partial series
for some weighted function $\mu(x),\ \ 0<\mu(x) \le1,
x\in[0,2\pi]$.

In [9] it is proved that for any given sequence of natural numbers
$\{\lambda_m \}_{m=1}^\infty$ with $\lambda_m \nearrow^\infty$
there exists a series by trigonometric system of the form
\begin{equation}
\sum_{k=1}^\infty C_k e^{ikx}, \ \  C_{-k}=\overline{C}_k,
\end{equation}
with
\[
  \left| \sum_{k=1}^m C_ke^{ikx} \right| \leq \lambda_m,\ \ x\in [0,2\pi],
  \ \ ,\ \ m=1,2,...,
\]
so that for each $\varepsilon>0$  a weighted function $\mu(x)$,
$$0<\mu(x) \le1, \left | \{ x\in[0,2\pi]: \mu(x)\not =1 \} \right |
<\varepsilon $$
can be constructed, so that the series (1.4) is
universal in the weighted space $L_\mu^1[0,2\pi]$ with respect
simultaneously to rearrangements as well as to subseries.

In this paper we prove the following results.
\begin{theorem}
 There exists a series of the form
\begin{equation}
\sum_{k=-\infty}^\infty{C_ke^{ikx} \ \  with \ \
\sum_{k=-\infty}^\infty \left | {C_k} \right|^q <\infty},\ \
\forall q>2
\end{equation}
such that for any number $\varepsilon>0$  a weighted function
$\mu(x)$, $ 0<\mu(x) \leq 1$, with
\begin{equation}
\left | \{ x\in[0,2\pi]:\mu(x)\not =1 \} \right | <\varepsilon
\end{equation}
can be constructed, so that the series $(1.5)$ is universal in
$L_{\mu}^1 [0,2\pi]$ with respect to rearrangements $.$
\end {theorem}
\begin{theorem}
There exists a series of the form $(1.5)$ such that for any number
$\epsilon>0$ a weighted function $\mu(x)$ with $(1.6)$ can be
constructed, so that the series $(1.5)$ is universal in
$L_\mu^1[0,2\pi]$ in the usual sense $.$
\end {theorem}

\par\par\bigskip
\par\par\bigskip

 \section{ BASIC LEMMA}
\par\par\bigskip
\par\par\bigskip

\begin{lemma}
 For any given numbers $0<\varepsilon <{\frac {1}{ 2}}$, $N_0>2$ and a step function
\begin{equation}
f(x)= \sum_{s=1}^q \gamma_s \cdot \chi_{\Delta_s} (x),
\end{equation}
where  $\Delta_s$ is an interval  of the form $\Delta_m^{(i)}=
\left[ {\frac {i-1}  {2^m}},{\frac {i} {2^m}} \right] $, $ 1\leq i
\leq 2^m $ and
\begin{equation}
|\gamma_s| \cdot \sqrt {|\Delta_s|}<\epsilon^3 \cdot \left( 8
\cdot \int_0^{2\pi} f^2(x)dx \right)^{-1},\ \ s=1,2,...,q.
\end{equation}
 there exists a measurable set $E \subset [0,2\pi]$ and a
polynomial $P(x)$ of the form
\[
P(x)= \sum_{N_0 \leq |k|<N} C_ke^{ikx}
\]
which satisfy the conditions:
\[
|E|> 2\pi- \varepsilon,
\]
\[
\int_E |P(x)-f(x)|dx<\varepsilon,
\]
\[
\sum_{N_0 \leq |k|<N} |C_k|^{2+\varepsilon}< \varepsilon,\ \
C_{-k}=\overline {C}_k
\]
\[
\max_{N_0 \leq m<N} \left[ \int_e \left | \sum_{N_0 \leq |k|\leq
m} C_k e^{ikx} \right | dx \right] <\varepsilon+\int_e |f_(x)|dx,
\]
for every measurable subset $e$  of $E$.

\end{lemma}
{\bf Proof } Let $0<\epsilon<{\frac {1}{ 2}}$ be an arbitrary
number.

Set
\begin{equation}
g(x)=1,\ \ if \ \ x \in [0,2\pi]\setminus \biggl[ {\frac
{\varepsilon \cdot \pi}{2}},  {\frac{3\varepsilon \cdot
\pi}{2}}\biggr] ;
\end{equation}
\[
 g(x)=1-{\frac{2} {\varepsilon}},\ \ if \ \ x \in \biggl[ {\frac { \varepsilon \cdot \pi}{ 2}},{\frac {3\varepsilon \cdot \pi}{ 2}}\biggr];
\]
We choose natural numbers $\nu_1$ and $N_1$ so large that the
following inequalities be satisfied:
\begin{equation}
{\frac {1} {2\pi}}\left|\int_0^{2\pi} g_1(t)e^{-ikt}dt \right|<
\frac {\varepsilon}{16 \cdot \sqrt{N_0}},\ \ |k|<N_0,
\end{equation}
where
\begin{equation}
g_1(x)=\gamma_1\cdot g(\nu_1\cdot x)\cdot \chi_{\Delta_1}(x).
\end{equation}
(By $\chi_E(x)$ we denote the characteristic function of the set
$E$.) We put
\begin{equation}
 E_1=\{ x\in \Delta_s:\ \ g_s(x)=\gamma_s \},
\end{equation}
By (2.3), (2.5) and (2.6) we have
\begin{equation}
|E_1|>2\pi\cdot (1-\epsilon)\cdot |\Delta_1|;\ \ g_1(x)=0,\ \
x\notin \Delta_1,
\end{equation}
\begin{equation}
\int_0^{2\pi}g_1^2(x)dx<{\frac {2}{ \epsilon}}\cdot
|\gamma_1|^2\cdot|\Delta_1|.
\end{equation}
Since the trigonometric system $\{e^{ikx}\}_{k=-\infty}^\infty$ is
complete in $L^2[0,2\pi]$, we can choose a natural number
$N_{1}>N_{0}$ so large that
\begin{equation}
\int_0^{2\pi}\left| \sum_{0\leq |k|<N_1}C_k^{(1)}e^{ikx}-
g_1(x)\right| dx \leq {\frac {\varepsilon} { 8}},
\end{equation}
where
\begin{equation}
C_k^{(1)}={\frac {1} {2\pi}}\int_0^{2\pi} g_1(t)e^{-ikt}dt .
\end{equation}
Hence by (2.4),(2.5) and (2.9) we obtain
\begin{equation}
\int_0^{2\pi}\left| \sum_{N_0\leq |k|<N_1}C_k^{(1)}e^{ikx}-
g_1(x)\right| dx\leq {\frac {\varepsilon} { 8}}+\left[ \sum_{0\leq
|k|<N_0} |C_k^{(1)}|^2 \right]^{\frac {1} { 2}}< {\frac
{\varepsilon}{4}};
\end{equation}
Now assume that the numbers $\nu_1<\nu_2<...\nu_{s-1}$,
$N_1<N_2<...<N_{s-1}$, functions $g_1(x),g_2(x),...,g_{s-1}(x)$
and the sets $E_1,E_2,....,E_{s-1}$ are defined. We take
sufficiently large natural numbers $\nu_s>\nu_{s-1}$ and
$N_s>N_{s-1}$ to satisfy
\begin{equation}
{\frac {1} {2\pi}}\left|\int_0^{2\pi} g_s(t)e^{-ikt}dt \right|<
\frac {\varepsilon} {16\cdot \sqrt{N_{s-1}}},\ \ \ \ 1\leq s \leq
q,\ \ \ \   |k|<N_{s-1},
\end{equation}
\begin{equation}
\int_0^{2\pi}\left| \sum_{0\leq |k|<N_s}C_k^{(s)}e^{ikx}-
g_s(x)\right| dx \leq {\frac {\varepsilon} { 4^{s+1}}},
\end{equation}
where
\begin{equation}
g_s(x)=\gamma_s\cdot g(\nu_s\cdot x)\cdot \chi_{\Delta_s}(x),\ \ \
\  C_k^{(s)}={\frac {1} {2\pi}}\int_0^{2\pi} g_s(t)e^{-ikt}dt .
\end{equation}
Set
\begin{equation}
E_s=\{ x\in \Delta_s:\ \ g_s(x)=\gamma_s \},
\end{equation}
Using the above arguments (see (2.16)-(2.18)), we conclude that
the function $g_s(x)$ and the set $E_s$ satisfy the conditions:
\begin{equation}
|E_s|>2\pi\cdot (1-\epsilon)\cdot |\Delta_s|;\ \ g_s(x)=0,\ \
x\notin \Delta_s,
\end{equation}
\begin{equation}
\int_0^{2\pi}g_s^2(x)dx<{\frac {2} {\epsilon}}\cdot
|\gamma_s|^2\cdot|\Delta_s|.
\end{equation}
\begin{equation}
\int_0^{2\pi}\left| \sum_{N_{s-1}\leq |k|<N_s}C_k^{(s)}e^{ikx}-
g_1(x)\right| dx < {\frac {\varepsilon} { 2^{s+1}}}.
\end{equation}
Thus, by induction we can define natural numbers
$\nu_1<\nu_2<...\nu_q$, $N_1<N_2<...<N_q$, functions
$g_1(x),g_2(x),...,g_q(x)$ and sets $E_1,E_2,....,E_q$ such that
conditions (2.14)- (2.16) are satisfied for all $s,\ \ 1\leq s
\leq q$. We define a set $E$ and a polynomial $P(x)$ as follows:
\begin{equation}
E=\bigcup_{s=1}^q E_s,
\end{equation}
\begin{equation}
P(x)=\sum_{N_0 \leq |k|< N} C_k e^{ikx}=\sum_{s=1}^q \left[
\sum_{N_{s-1}\leq |k|<N_s} C_k^{(s)}e^{ikx} \right],\ \
C_{-k}=\overline {C}_k,
\end{equation}
where
\begin{equation}
 C_k=C_k^{(s)} \ \  for\ \ N_{s-1}\leq |k|<N_s,\ \ s=1,2,...,q, \ \  N=N_q-1.
\end{equation}
By Bessel's inequality and (2.3), (2.14) for all $s\in [1,q]$ we
get
\begin{equation}
\left[ \sum_{N_{s-1}\leq |k|<N_s} |C_k^{(s)}|^2 \right]^{\frac
{1}{2}} \leq \left[ \int_o^{2\pi}g_s^2(x)dx \right]^{\frac
{1}{2}}\leq
\end{equation}
\[
\leq {\frac {2} { \sqrt{\varepsilon}}}\cdot |\gamma_s| \cdot \sqrt
{|\Delta_s|},\ \ s=1,2,...,q.
\]
From (2.3), (2.12) and (2.13) it follows that
\begin{equation}
|E|> 2\pi- \varepsilon.
\end{equation}
Taking relations (2.1), (2.3), (2.10), (2.12), (2.18) - (2.21) we
obtain
\begin{equation}
\int_E |P(x)-f(x)|dx \leq \sum_{s=1}^q \left[ \int_E \left|
\sum_{N_{s-1}\leq |k|<N_s} C_k^{(s)}e^{ikx}-g_s(x) \right|
dx\right]<\varepsilon
\end{equation}
By (2.1), (2.2), (2.20) -(2.21)for any $k\in [N_0,N]$ we have
\begin{equation}
\sum_{N_0\leq |k|< N} |C_k|^{2+ \epsilon} \leq \max_{N_0 \leq k
\leq N} |C_k|^\epsilon \cdot \sum_{k=N_0}^N |C_k|^2 \leq
\end{equation}
\[
\leq \max_{1 \leq s \leq q}\left[\sqrt {\frac {8} {\epsilon}}
\cdot |\gamma_s| \cdot \sqrt {|\Delta_s|}\right] \cdot
\sum_{s=1}^q \left[ \sum_{N_{s-1}\leq |k|<N_s} |C_k^{(s)}|^2
\right] \leq
\]
\[
 \leq \max_{1 \leq s \leq q}\left[ \sqrt {\frac {8} {\epsilon}}  \cdot |\gamma_s| \cdot \sqrt {|\Delta_s|}\right] \cdot {\frac {8} { \epsilon}} \cdot \sum_{s=1}^q |\gamma_s|^2 \cdot|\Delta_s| \leq
\]
\[
 \leq\max_{1 \leq s \leq q}\left[ \sqrt  {\frac {8} {\epsilon}} \cdot |\gamma_s| \cdot \sqrt {|\Delta_s|}\right] \cdot {\frac {8} {\epsilon}} \cdot \left[ \int_0^1 f^2(x)dx \right]< \epsilon;
\]
That is, the statements 1) - 3) of Lemma are satisfied. Now we
will check the fulfillment of statement 4) of Lemma. Let $N_0 \leq
m<N$, then for some $s_0,\ \ 1 \leq s_0 \leq q, \ \ \left( N_{s_0}
\leq m< N_{s_0+1} \right) $ we will have (see (2.20) and (2.21))
\begin{equation}
\sum_{N_0 \leq |k|\leq m} C_k e^{ikx}=\sum_{s=1}^{s_0} \left[
\sum_{N_{s-1}\leq |k|<N_s} C_k^{(s)}e^{ikx}
\right]+\sum_{N_{s_0-1}\leq |k|\leq m} C_k^{(s_0+1)} e^{ikx}.
\end{equation}
Hence and from (2.1), (2.2), (2.3), (2.18), (2.19) and (2.22) for
any measurable set $e \subset E$ we obtain
\[
\int_e \left | \sum_{N_{s-1}\leq |k|\leq m} C_k e^{ikx} \right |
dx \leq
\]
\[
\leq \sum_{s=1}^{s_0} \left[\int_e \left | \sum_{N_{s-1}\leq
|k|<N_s} C_k^{(s)}e^{ikx}-g_s(x) \right|dx \right] +
\]
\[
+\sum_{s=1}^{s_0} \int_e |g_s(x) | dx +\int_e \left|
\sum_{N_{s_0-1}\leq |k|\leq m} C_k^{(s_0+1)} e^{ikx} \right|dx<
\]
\[
<\sum_{s=1}^{s_0} {\frac {\varepsilon} {2^{s+1}}}+\int_e
|f(x)|dx+{\frac {2} {\sqrt {\varepsilon}}}  \cdot |\gamma_{s_0+1}|
\cdot \sqrt {|\Delta_{s_0+1}|} <
\]
\[
< \int_e |f(x)|dx+ \varepsilon.
\]
\par\par\bigskip
\par\par\bigskip

\section{PROOF OF THEOREMS }
\par\par\bigskip
\par\par\bigskip

{\bf Proof of Theorem 1.4 }
 Let
\begin{equation}
f_1(x), f_2(x),...,f_n(x), \ \  x \in [0,2\pi]
\end{equation}
be a sequence of all step functions, values and constancy interval
endpoints of which are rational numbers. Applying Lemma
consecutively, we can find a sequence $\{ E_s\}_{s=1}^\infty $ of
sets and a sequence of polynomials
\begin{equation}
P_s(x)=\sum_{N_{s-1}\leq |k|<N_s} C_k^{(s)}e^{ikx}
\end{equation}
\[
 1=N_0<N_1<...<N_s<....,\ \ s=1,2,....,
\]
which satisfy the conditions:
\begin{equation}
| E_s| >1-2^{-2(s+1)} ,\ \  E_s \subset [0,2\pi],
\end{equation}
\begin{equation}
\int_{E_s}|P_s(x)-f_s(x)|dx<2^{-2(s+1)},
\end{equation}
\begin{equation}
 \sum_{N_{s-1}\leq |k|<N_s}\left |C_k^{(s)}\right|^{2+2^{-2s}}< 2^{-2s},\ \ C_{-k}^{(s)}=\overline {C}_k^{(s)}
\end{equation}
\begin{equation}
\max_{N_{s-1}\leq p<{N_s}} \left[ \int_e \left | \sum_{N_{s-1}\leq
|k|\leq p}C_k e^{ikx} \right | dx \right] <2^{-2(s+1)}+\int_e
|f_s(x)|dx,
\end{equation}

for every measurable subset $e$  of $E_s$.

Denote
\begin{equation}
\sum_{k=-\infty}^\infty C_k e^{ikx}=\sum_{s=1}^\infty \left[
\sum_{N_{s-1}\leq |k|<N_s}  C_k^{(s)}e^{ikx} \right],
\end{equation}
where $ C_k=C_k^{(s)}$ for $N_{s-1}\leq |k|<N_s$, $s=1,2,...$.

Let  $\varepsilon$ be an arbitrary positive number. Setting
\[
\Omega_n = \bigcap_{s=n}^\infty E_s,\ \  n=1,2,....;
\]
\begin{equation}
 E=\Omega_{n_0} = \bigcap_{s=n_0}^\infty E_s,\ \   n_0=[\log_{1/2} \varepsilon]+1;
\end{equation}
\[
B=\bigcup _{n=n_0} ^\infty \Omega_n =\Omega_{n_0} \bigcup \left(
\bigcup _{n=n_0+1}^ \infty  \Omega_n \setminus \Omega_{n-1}
\right).
\]

It is clear (see (3.3)) that $| B |=2\pi$  and $| E | >2\pi-
\varepsilon$.

We define a function $\mu(x)$ in the following way:
\begin{equation}
\mu(x)=1\ \ for \ \ x \in E \cup ([0,2\pi] \setminus B);
\end{equation}
\[
\mu(x)=  \mu_{ n} \ \  for  \ \ x \in \Omega_{n} \setminus
\Omega_{n-1},\ \  n\geq n_0+1,
\]
where
\begin{equation}
 \mu_n=\left[ 2^{4n}\cdot \prod_{s=1}^n h_s \right]^{-1};
\end{equation}
\[
  h_s=|| f_s(x)||_C+ \max_{N_{s-1}\leq p<{N_s}} \Vert  \sum_{N_{s-1}\leq |k|\leq p}  C_k^{(s)}e^{ikx} \Vert_{C}+1,
\]

where
\[
 ||g(x)||_C=\max_{x\in [0,2\pi]} |g(x)|,
\]
$g(x)$ is a continuous function on  $[0,2\pi]$.

From (3.5),(3.7)-(3.10) we obtain

(A) --   $ 0<\mu(x) \le1, \mu(x)$ is a measurable function and
\[
\left | \{x\in[0,2\pi]:\mu(x)\not =1\} \right|<\varepsilon.
\]

(B) --   $\displaystyle \sum_{k=1}^\infty |C_k|^q<\infty,\ \
\forall q>2.$

Hence, obviously we have
\begin{equation}
 \lim_{k\to \infty}{C_k}=0.
\end{equation}
It follows  from (3.8)-(3.10) that for all $s \geq n_0$  and $p
\in \left[ N_{s-1},N_s \right)$
\begin{equation}
\int_{[0,2\pi] \setminus \Omega_s} \left|  \sum_{N_{s-1}\leq
|k|\leq p}  C_k^{(s)}e^{ikx} \right| \mu(x) dx=
\end{equation}
\[
=\sum_{n=s+1}^ \infty \left[\int_{\Omega_n \setminus \Omega_{n-1}}
\left|  \sum_{N_{s-1}\leq |k|\leq p}  C_k^{(s)}e^{ikx}\right|
\mu_n dx \right] \leq
\]
\[
\leq \sum_{n=s+1}^ \infty2^{-4n} \left[\int_0^ {2\pi} \left|
\sum_{N_{s-1}\leq |k|\leq p}  C_k^{(s)}e^{ikx}\right| h_s^{-1} dx
\right]<2^{-4s}.
\]
By (3.4), (3.8)-(3.10) for all  $s \geq n_0$ we have
\begin{equation}
\int_0^{2\pi} \left| P_s(x)-f_s(x) \right|\mu(x)dx=\int_{\Omega_s}
\left| P_s(x)-f_s(x) \right|\mu(x)dx+
\end{equation}
\[
 +\int_{[0,2\pi] \setminus \Omega_{s}} \left| P_s(x)-f_s(x) \right|\mu(x)dx =2^{-2(s+1)}+
\]
\[
+\sum_{n=s+1}^\infty \left[\int_{\Omega_n \setminus \Omega_{n-1}}
\left| P_s(x)-f_s(x)  \right| \mu_n dx\right] \leq 2^{-2(s+1)}+
\]
\[
+ \sum_{n=s+1}^ \infty 2^{-4s}\left[ \int_0^ {2\pi} \left(\left|
f_s(x) \right| +\left| \sum_{N_{s-1}\leq |k|<N_s}
C_k^{(s)}e^{ikx} \right| \right) h_s^{-1}dx \right ]<
\]
\[
<2^{-2(s+1)}+2^{-4s}<2^{-2s}.
\]
Taking relations (3.6), (3.8)- (3.10) and (3.12) into account we
obtain that for all $p \in \left[ N_{s-1},N_s \right)$ and $s \geq
n_0+1$
\begin{equation}
\int_0^{2\pi} \left| \sum_{N_{s-1}\leq |k|\leq p}
C_k^{(s)}e^{ikx} \right| \mu(x) dx=
\end{equation}
\[
=\int_{\Omega_s} \left| \sum_{N_{s-1}\leq |k|\leq p}
C_k^{(s)}e^{ikx} \right| \mu(x) dx+
\]
\[
+\int_{[0,2\pi] \setminus \Omega_s} \left| \sum_{N_{s-1}\leq
|k|\leq p}  C_k^{(s)}e^{ikx} \right| \mu(x) dx<
\]
\[
< \sum_{n=n_0+1}^ s \left[\int_{\Omega_n \setminus \Omega_{n-1}}
\left| \sum_{N_{s-1}\leq |k|\leq p}  C_k^{(s)}e^{ikx} \right| dx
\right]\cdot \mu_n+2^{-4s}<
\]
\[
<\sum_{n=n_0+1}^ s \left( 2^{-2(s+1)}+\int_{\Omega_n \setminus
\Omega_{n-1}} |f_s(x)|dx \right) \mu_n +2^{-4s}  =
\]
\[
=2^{-2(s+1)} \cdot \sum_{n=n_0+1}^ s \mu_n+\int_{\Omega_s}
|f_s(x)|\mu(x)dx +2^{-4s}<
\]
\[
<\int_0^{2\pi} |f_s(x)|\mu(x)dx +2^{-4s}.
\]
Let $ f(x) \in L_{\mu}^1 [0,2\pi]$ , i. e.$ \int_0^{2\pi} |f(x)|
\mu(x) dx<\infty$ .

It is easy to see that we can choose a function $f_{\nu_1}(x)$
from the sequence (3.1) such that
\begin{equation}
\int_0^{2\pi} \left| f(x)- f_{\nu_1}(x) \right|\mu(x)dx<2^{-2},\ \
\nu_1 > n_0+1.
\end{equation}
Hence, we have
\begin{equation}
\int_0^{2\pi} \left| f_{\nu_1}(x)
\right|\mu(x)dx<2^{-2}+\int_0^{2\pi} |f(x)|\mu(x)dx.
\end{equation}
From (2.1), (A), (3.13) and (3.15)  we obtain with $m_1=1$
\begin{equation}
\int_0^{2\pi} \left| f(x)- \left [ P_{\nu_1}(x)+C_{m_1}e^{im_1x}
\right] \right|\mu(x)dx \leq
\end{equation}
\[
\leq \int_0^{2\pi} \left| f(x)- f_{\nu_1}(x) \right|\mu(x)dx+
\]
\[
+\int_0^{2\pi} \left| f_{\nu_1}(x)-P_{\nu_1}(x) \right|\mu(x)dx+
\]
\[
+\int_0^{2\pi} \left| C_{m_1}e^{im_1x}\right|\mu(x)dx <2\cdot
2^{-2}+2\pi\cdot\left| C_{m_1}\right|.
\]
Assume that numbers
$\nu_1<\nu_2<...<\nu_{q-1};m_1<m_2<...<m_{q-1}$ are chosen in such
a way that the following condition is satisfied:
\begin{equation}
\int_0^{2\pi} \left| f(x)- \sum_{s=1}^j \left [
P_{\nu_s}(x)+C_{m_s}e^{im_sx} \right] \right|\mu(x)dx<
\end{equation}
\[
<2\cdot 2^{-2j}+2\pi \cdot\left|C_{m_j}\right|, \ \ 1\leq j \leq
q-1 .
\]
We choose a function $f_{\nu_q}(x)$ from the sequence (3.1) such
that
\begin{equation}
\int_0^{2\pi} \left| \left( f(x)- \sum_{s=1}^{q-1} \left [
P_{\nu_s}(x)+C_{m_s}e^{im_sx} \right]  \right)-f_{n_q}(x)\right|
\mu(x)dx< 2^{-2q},
\end{equation}
where $ \nu_q>\nu_{q-1};\ \  \nu_q>m_{q-1}$

This with (3.18) imply
\begin{equation}
\int_0^{2\pi} \left| f_{\nu_q}(x) \right| \mu(x)dx<2^{-2q}+2\cdot
2^{-2(q-1)}+2\pi\cdot\left|C_{m_{q-1}} \right| =
\end{equation}
\[
=9 \cdot 2^{-2q}+ 2\pi\cdot\left|C_{m_{q-1}} \right|.
\]

By (3.13), (3.14) and (3.20) we obtain
\begin{equation}
\int_0^{2\pi} \left| f_{\nu_q}(x)- P_{\nu_q}(x)
\right|\mu(x)dx<2^{-2\nu_q},
\end{equation}
\[
\ \ P_{\nu_q}(x)=\sum_{N_{\nu_q-1}\leq |k|<N_{\nu_q}}
C_k^{(\nu_q)}e^{ikx}.
\]
\begin{equation}
\max_{N_{\nu_q-1} \leq p<N{\nu_q}}  \int_0^{2\pi} \left|
\sum_{k=N_{\nu_q-1}}^ p C_k^{(\nu_q)}e^{ikx} \right | \mu(x) dx<10
\cdot 2^{-2q}+2\pi\cdot\left | C_{m_{q-1}} \right |.
\end{equation}
Denote
\begin{equation}
m_q= \min \left\{ n \in N: n \notin \left\{  \left\{  \{ k
\}_{k=N_{\nu_s-1}}^{N_{\nu_s}-1} \right\}_{s=1}^q \cup \{
m_s\}_{s=1}^{q-1}\right\} \right\}.
\end{equation}
From (2.1), (A), (3.19) and (3.21) we have
\begin{equation}
\int_0^{2\pi} \left|  f(x)- \sum_{s=1}^q \left [
P_{\nu_s}(x)+C_{m_s}e^{im_sx} \right] \right| \mu(x)dx\leq
\end{equation}
\[
\leq \int_0^{2\pi} \left| \left( f(x)- \sum_{s=1}^{q-1} \left [
P_{\nu_s}(x)+C_{m_s}e^{im_sx} \right]  \right)-f_{\nu_q}(x)\right|
\mu(x)dx+
\]
\[
+\int_0^{2\pi} \left| f_{\nu_q}(x)-P_{\nu_q}(x)\right| \mu(x)dx+
\]
\[
+\int_0^{2\pi} \left| C_{m_q}e^{im_qx}\right| \mu(x)dx<2 \cdot
2^{-2q}+2\pi\cdot\left| C_{m_q} \right|.
\]
Thus, by induction we on $q$ can choose from series (3.7) a
sequence of members
\[
C_{m_q}e^{im_qx} ,\ \  q=1,2,...,
\]
and a sequence of polynomials
\begin{equation}
P_{\nu_q}(x)=\sum_{N_{\nu_q-1}\leq |k|<N_{\nu_q}}
C_k^{(\nu_q)}e^{ikx},\ \  N_{n_q-1}>N_{n_{q-1}},\ \  q=1,2,....
\end{equation}
such that conditions (3.22) - (3.24) are satisfied for all $q\geq
1.$

Taking account the choice of $P_{\nu_q}(x)$ and $C_{m_q}e^{im_qx}$
(see (3.23) and (3.25)) we conclude that the series
\[
\sum_{q=1}^\infty \left[ \sum_{N_{\nu_q-1}\leq |k|<N_{\nu_q}}
C_k^{(\nu_q)}e^{ikx}+C_{m_q}e^{iqx} \right ]
\]
is obtained from the series (3.7) by rearrangement of members.
Denote this series by $\sum C_{\sigma(k)}e^{i\sigma(k) x}.$

It follows from (3.11), (3.22) and (3.24) that the series $\sum
C_{\sigma(k)}e^{i \sigma(k)x}$ converges to the function $f(x)$ in
the metric $L_{\mu}^1[0,2\pi]$, i.e. the series (3.7) is universal
with respect to rearrangements (see Definition 1.1).

{\bf The Theorem 1.4 is proved.}
\[
\]
{\bf Proof of the Theorem 1.5}

Applying Lemma  consecutively, we can find a sequence $\{
E_s\}_{s=1}^\infty $ of sets and a sequence of polynomials
\begin{equation}
P_s(x)=\sum_{N_{s-1}\leq |k|<N_s} C_k^{(s)}e^{ikx} ,\ \
C_{-k}^{(s)}=\overline {C}_k^{(s)}
\end{equation}
\[
 1=N_0<N_1<...<N_s<....,\ \ s=1,2,....,
\]

which satisfy the conditions:
\begin{equation}
\left| E_s\right| >1-2^{-2(s+1)} ,\ \  E_s\subset [0,2\pi],
\end{equation}
\begin{equation}
 \sum_{N_{s-1}\leq |k|<N_s}\left |C_k^{(s)}\right|^{2+2^{-2s}}< 2^{-2s},
\end{equation}
\begin{equation}
  \int_{E_n} \left | f_n(x)-\sum_{s=1}^n P_s(x)\right | dx <2^{-n},\ \ n=1,2,...,
\end{equation}
where $\{ f_n(x)\}_{ n=1}^\infty,\ \  x\in [0,2\pi]$ be a sequence
of all step functions, values and constancy interval endpoints of
which are rational numbers.

Denote
\begin{equation}
\sum_{k=-\infty}^\infty C_k e^{ikx}=\sum_{s=1}^\infty \left[
\sum_{N_{s-1}\leq |k|<N_s}  C_k^{(s)}e^{ikx} \right],
\end{equation}
where $ C_k=C_k^{(s)}$ for $N_{s-1}\leq |k|<N_s$, $s=1,2,...$.

It is clear (see (3.28)) that
\[
\displaystyle \sum_{k=1}^\infty |C_k|^q<\infty, \ \ \forall q>2.
\]

Repeating reasoning of Theorem 1 a weighted function  $\mu (x),\ \
0<\mu(x)\leq 1$ can constructed so that the following condition is
satisfied:
\begin{equation}
\int_0^{2\pi} \left| f_n(x)-\sum_{s=1}^n P_s(x) \right|\cdot \mu
(x) dx< 2^{-2n},\ \ n=1,2,...
\end{equation}
For any function $f(x) \in L_{\mu}^1[0,1]$ we can choose a
subsystem $\{f_{n_\nu}(x)\}_{\nu=1}^\infty $ from the sequence
(3.1) such that
\begin{equation}
\int_0^{2\pi} \left| f(x)- f_{n_\nu}(x) \right|\mu(x)dx<2^{-2\nu}.
\end{equation}
From (3.30)-(3.32) we conclud
\[
\int_0^{2\pi} \left|  f(x)- \sum_{|k|\leq M_\nu}C_ke^{ikx}
\right| \mu(x)dx =
\]
\[
\int_0^{2\pi} \left|  f(x)- \sum_{s=1}^{n_\nu} \left[
\sum_{N_{s-1}\leq |k|<N_s}  C_k^{(s)}e^{ikx} \right] \right|
\mu(x)dx \leq
\]
\[
\leq \int_0^{2\pi} \left| f(x)- f_{\nu_k}(x) \right| \cdot
\mu(x)dx+
\]
\[
+\int_0^{2\pi} \left| f_{\nu_k}(x)-\sum_{s=1}^{\nu_k} P_s(x)
\right|\cdot \mu (x) dx<2^{-2k}+2^{-2{\nu_k}}
\]
where $M_{\nu}=N_{n_{\nu}}-1$.

Thus, the series (3.30) is universal in $L_\mu^1[0,1]$ in the
sense of usual (see Definition 1.2).

{\bf The Theorem 1.5 is proved.}
\par\par\bigskip

The author thanks Professor M.G.Grigorian for his attention to
this paper.

\newpage

\bibliographystyle{amsplain}
\par\par\bigskip
\par\par
\end{document}